\documentclass[11pt]{article}

\usepackage{amsthm,amssymb}

\newtheorem{thm}{Theorem}[section]
\newtheorem{prop}{Proposition}[section]
\newtheorem{cor}{Corollary}[section]
\newtheorem{rem}{Remark}[section]
\newtheorem{lem}{Lemma}[section]

\newcommand{\dist}{{\mathrm{dist}\,}}
\newcommand{\eps}{\varepsilon}

\newcommand{\vol}{{\mathrm{Vol}\,}}

\newcommand{\rank}{{\mathrm{rank}}}

\newcommand{\const}{\mathrm{const}}
\newcommand{\R}{{\mathbb{R}}}

\newcommand{\mC}{{\mathbb{C}}}
\newcommand{\Z}{{\mathbb{Z}}}
\newcommand{\T}{{\mathbb{T}}}
\newcommand{\N}{{\mathbb{N}}}

\newcommand{\calB}{\mathcal{B}}

\newcommand{\Real}{{\mathrm{Re}\,}}

\title{ Topological approach to the generalized $n$-center problem 
\author{Sergey Bolotin\footnote{Moscow Steklov Mathematical Institute and University of Wisconsin-Madison.} \ 
and Valery Kozlov\footnote{Moscow Steklov Mathematical Institute.}} }

\date{}

\begin{document}

\maketitle

  \begin{abstract}
We consider  a  natural Hamiltonian system with  two degrees of freedom and Hamiltonian $H=\|p\|^2/2+V(q)$.
The configuration space $M$ is a closed surface (for noncompact $M$ certain conditions at infinity are required).
It is well known that if the potential energy $V$  has $n>2\chi(M)$ Newtonian singularities,
then the system is not integrable and has positive topological entropy on energy levels
$H=h>\sup V$. We generalize this result to the case when the potential energy has several singular points $a_j$
of type $V(q)\sim -\dist(q,a_j)^{-\alpha_j}$.
Let $A_k=2-2k^{-1}$, $k\in\N$, and let  $n_k$ be the number of singular points
with $A_k\le \alpha_j<A_{k+1}$.
We prove that if
$$
\sum_{2\le k\le\infty}n_kA_k>2\chi(M),
$$
then the system has  a compact chaotic invariant set
of noncollision trajectories on any energy level $H=h>\sup V$.
This result is purely topological: no analytical properties of the potential,
except the presence of singularities, are involved. The proofs are based
on the generalized Levi-Civita regularization and elementary topology of coverings.
As an example, the plane $n$ center problem is considered.
\end{abstract}

\section{Introduction}

\label{sec:intro}

Let $M$ be a  connected 2-dimensional manifold -- the configuration space of a  natural Hamiltonian
system with two degrees of freedom. Passing to a 2-sheet covering we may assume that $M$ is oriented.
The Hamiltonian $H$ on $T^*M$ is quadratic
in the momentum:
\begin{equation}
\label{eq:H}
H(q,p)=\frac12\|p\|^2+V(q),\qquad p\in T_q^*M,
\end{equation}
where $\|\cdot \|$ is a Riemannian metric on $M$. The corresponding Lagrangian
is\footnote{We use the same notation for the norm of a vector and a covector.}
\begin{equation}
\label{eq:L}
L(q,\dot q)=\frac12\|\dot q\|^2-V(q).
\end{equation}
Trajectories of the system satisfy the Newton equation
\begin{equation}
\label{eq:Newton}
\frac{D\dot q}{dt}=-\nabla V(q),
\end{equation}
where $D/dt$ is the covariant derivative and $\nabla V$ the gradient vector.

We  assume that the metric  is smooth\footnote{Smooth means of class at least $C^3$.} on $M$,
and  the potential energy $V$ is smooth except for a  finite set of singular points $\Delta=\{a_1,\dots,a_n\}\subset  M$.
More precisely, $V$
is smooth on $M\setminus \Delta$ and in small balls
\begin{equation}
 \label{eq:ball}
B_j=B(a_j,\eps)=\{q\in M:d(x,a_j)\le \eps\}
 \end{equation}
  it has the form
\begin{equation}
V(q)=-\frac{f_j(q)}{d (q,a_j)^{\alpha_j}}+U_j(q),\qquad \alpha_j>0,\quad m_j=f_j(a_j)> 0,
\label{eq:V}
\end{equation}
where the functions  $f_j$ and $ U_j $ are smooth   on $ B_j$.
The distance $ d (q, a_j) $ is measured in the Riemannian metric $\|\cdot\|$.
The configuration space of the system is $\hat M=M\setminus\Delta$, and the phase space is $T^*\hat M$.
 
The orders of singularities $\alpha_j$ are arbitrary positive numbers.
The most physical are   Newtonian singularities with $\alpha_j=1$. Singularities with $\alpha_j\ge 2$
are called strong force singularities  \cite{Gordon}.
As discovered already by Poincar\'e, for  strong  force singularities it is easy
to prove the existence of periodic  solutions by variational methods, see   \cite{Ambr,Gordon}.
The reason is that trajectories colliding with such singularities have infinite action.
The case of  singularities with $0<\alpha_j<2$ is more difficult. We call a singularity weak if
$0<\alpha_j<1$ and moderate if $1<\alpha_j<2$. Newtonian singularities with $\alpha_j=1$
and  singularities with $\alpha_j=2$ are critical. We call singularities with $\alpha_j=2$ Jacobi singularities.
Jacobi studied the $n$-body space problem with the potential of degree $-2$ 
and discovered simple behavior of the moment of inertia of the system. When $\alpha_j>2$ we say that the singularity is strong.
For proving the existence of chaotic trajectories there is a considerable difference between
strong and Jacobi singularities. Thus
$$
\Delta=\Delta_{weak}\cup\Delta_{newt}\cup\Delta_{mod}\cup\Delta_{jac}\cup \Delta_{strong}.
$$

A standard example is the generalized $n$ center problem in $\R^2$:
\begin{equation}
\label{eq:ncent}
H(q,p)=\frac12|p|^2+V(q),\qquad V(q)=-\sum_{j=1}^n\frac{m_j}{|q-a_j|^{\alpha_j}}+U(q),\qquad   m_j>0,
\end{equation}
where $U$ is smooth and bounded above on $\R^2$, and the metric is Euclidean. Usually it is assumed that the orders $\alpha_j$
are equal, but we do not impose this assumption.

When all singularities are Newtonian ($\alpha_j=1$) and $n\ge 3$  the $n$-center problem has  a
chaotic invariant set on any energy level $H=h>\sup V$  for purely topological reasons,
see  \cite{Bol:ncent,Knauf,Bol:geo}.
For $n=2$ the classical 2-center problem with $\alpha_j=1$ and $U=0$ has a quadratic in momentum first integral.
For strong   singularities  ($\alpha_j>2$) the $n$ center problem has a chaotic invariant set on the level $H=h>\sup V$
for $n\ge 2$, again for purely topological reasons.
In the present paper we prove topological sufficient conditions for chaotic behavior of
the $n$ center problem for arbitrary $\alpha_j>0$, generalizing some results of \cite{Castel,Terra1,Terra2}.
For recent references see \cite{Castel}.

We consider trajectories of the system (\ref{eq:H}) with a fixed value of the total energy
\begin{equation}
\label{eq:h}
H=h>\max_M V.
\end{equation}
When $ h <\sup_M V $, motion of the system occurs in the domain of possible motion $D_h=\{q\in M:V\le h\}$.
Then our main theorem does not hold. In particular, they do not hold for repelling singularities, then $\sup V=+\infty$.
Then only much weaker results can be proved, see Remark \ref{rem:libr}.

The problem of integrability of system (\ref{eq:H}) on the energy level (\ref{eq:h}) was discussed in \cite{Bol-Koz:Izv}.
One of the results of \cite{Bol-Koz:Izv} is the following theorem. Let $\chi(M)$ be the Euler characteristics of $M$.

\begin{thm}\label{thm:main}
Let $M$ be a closed manifold. Suppose there are only moderate and Newtonian singularities with  $1\le \alpha_j<2$. If
\begin{equation}
\label{eq:sum}
\sum_{j=1}^n\alpha_j>2\chi(M),
\end{equation}
then  the Hamiltonian system has no nonconstant polynomial in momenta and smooth
in coordinates conditional first integrals (Birkhoff conditional integrals)  on the energy level
$\{H=h\}\subset T^*\hat M$ for $h>\max_M V$.
\end{thm}

If there are no singularities ($\Delta=\emptyset$), this is the result of   \cite{Kozlov:top}:
a natural analytic system on a surface with genus greater than one can not  have nonconstant analytic conditional integrals.
For smooth systems, the same argument implies nonexistence of Birkhoff conditional integrals.
See  also \cite{Kolokol} for a direct proof. 
For Newtonian singularities condition (\ref{eq:sum}) gives $n=\#\Delta>2\chi(M)$.
In this case Theorem \ref{thm:main} was proved in \cite{Bol:sing}.

Theorem \ref{thm:main} holds also for systems with gyroscopic (or magnetic)
forces when the symplectic form is
\begin{equation}
\label{eq:Omega}
dp\wedge dq+\pi^*\omega,
\end{equation}
where $\omega$ is a closed 2-form on $M$ and $\pi:T^*M\to M$ is the projection,
see \cite{Bol-Koz:Izv}.
However, if the form $\omega$ is nonexact then under condition (\ref{eq:sum})
the system has no nonconstant Birkhoff integrals but it may have  nonconstant analytic
in momentum conditional integrals on an energy level (\ref{eq:h})
and hence no chaotic trajectories, see e.g.\ \cite{Bol:gyro}.

For example, suppose there are  no singularities  and $\chi(M)<0$. For a given metric and a constant  $h$ 
choose the potential energy so that the  Jacobi metric (\ref{eq:Jac}) has constant negative curvature.  Set
$\omega=c\Omega$, where  $\Omega$ is the area form of the Jacobi metric.
The for a large constant $c$  all trajectories with energy $h$ are periodic,
so the system has an analytic conditional integral on $\{H=h\}$, but no nonconstant Birkhoff integrals.

For natural systems (\ref{eq:H}) we suspect that  condition (\ref{eq:sum}) implies the existence
of chaotic  trajectories, but in general this is not proved.

We will give topological  conditions slightly stronger than (\ref{eq:sum})
which imply  the existence of chaotic trajectories,
in particular positiveness of the topological entropy. A simple  corollary of the main theorem is the following:

\begin{thm}\label{thm:simple}
Let $M$ be a closed manifold. If $1\le \alpha_j<2$ and
\begin{equation}
\label{eq:n>}
n>2\chi(M),
\end{equation}
then for $h>\max_M V$ the  system has a compact invariant set
with positive topological entropy on the energy level
$\{H=h\}\subset T^*\hat M$.
\end{thm}

When all singularities are Newtonian, Theorem {thm:simple} follows from a result in   \cite{Bol:sing}.
The proof is based on  the global Levi-Civita regularization
\cite{L-C}. Moderate singularities with $1<\alpha_j<2$ are mostly nonregularizable.
However, the method of the proof of Theorem \ref{thm:simple} is essentially   the same as in \cite{Bol:sing}.
We will prove:

\begin{prop}\label{prop:closed}
There is a closed surface $N$, a smooth $K$-sheet covering $\phi:N\to M$,
branched over $\Delta$,
and a smooth Riemannian metric on $N$ such that:
\begin{itemize}
\item
The covering has $K=2$ sheets if $n$ is even and $K=4$ sheets if $n$ is odd.
\item
The Euler characteristics of $N$ is
$$
\chi(N)=K\Big(\chi(M)-\frac n2\Big)<0.
$$
\item
The map $\phi$
takes minimal geodesics on $N$ to trajectories
of the system on the energy level $\{H=h\}$ (with changed parametrization) not colliding with $\Delta$.
An exception are regularizable Newtonian singularities:  some minimal geodesics may pass through
$ \phi^{-1}(\Delta_{newt})$ and correspond to trajectories on $M$ reflecting from  a Newtonian singularity.
\end{itemize}
\end{prop}

A geodesic on $N$ is called minimal if it minimizes the distance between any two points on its
lift to the universal covering of $N$. Note that nonminimal geodesics on $N$ may not correspond
to trajectories of energy $h$ on $M$.

In particular, any nontrivial homotopy class of closed curves
in $N$ contains a minimal geodesic which is projected to a periodic orbit with energy $h$
having  no collisions with $\Delta_{mod}$.

The surface $N$ is a sphere with more than one handle. For example, if $M=S^2$ and $n=2k$, the
Euler characteristics of $N$ is $\chi(N)=2(2-k)$
and the genus is  $g=1-\chi/2=k-1$,  so $N$ is a sphere with $k-1$ handles.

To deduce Theorem \ref{thm:simple} from Proposition \ref{prop:closed} we recall that the geodesic flow on a closed surface $N$
with genus greater  than one has  a compact chaotic invariant set.
This was essentially known to Morse and Hedlund \cite{Morse,Hedlund} who proved the existence
of an infinite number of minimal heteroclinic geodesics joining closed geodesics. One can show \cite{Bol:geo}
that these heteroclinics are topologically transverse. A rigorous proof of positiveness of the   topological entropy
was given by Dinaburg \cite{Dinaburg}. In fact also the set of minimal geodesics forms a compact
chaotic invariant set with positive topological entropy, see   \cite{Katok,KH}.
As mentioned above, the corresponding trajectories in $M$
may have regularizable collisions with $\Delta_{newt}$, but then there is another set  of positive topological entropy
and no collisions.

We  briefly  recall what happens in the presence of regularizable Newtonian singularities, see also \cite{Bol:sing}.
There is a sheet exchanging involution $\sigma:N\to N$,    $\phi\circ\sigma=\phi$, preserving the metric,
such that the set of fixed
points of $\sigma$ is  $\phi^{-1}(\Delta)$.
 The global Levi-Civita regularization $\phi:N\to M$
completely removes  Newtonian singularities so  trajectories  $\gamma$ in   $N$
may freely pass  through $\phi^{-1}(\Delta_{newt})$.
If $\gamma(0)\in \phi^{-1}(\Delta_{newt})$, the   minimizer will be $\sigma$-reversible:
$\gamma(-t)=\sigma\gamma(t)$.
The corresponding trajectory $q(t)=\phi(\gamma(t))$ in $M$ will have a reflection from  a Newtonian singularity:
$q(0)\in\Delta_{newt}$ and
$q(-t)=q(t)$. However, only a small portion of homotopy classes
of closed curves  in $N$ are $\sigma$-reversible. A minimizer in a nonreversible homotopy class will be projected
to a noncollision trajectory of energy $h$.
Similarly, a big part of   minimal (on the universal  covering of $N$) nonperiodic
geodesics will have no collisions with Newtonian singularities.

 In the next section we formulate a generalization of Theorem \ref{thm:simple}.

\section{Main result}

Keeping in mind the $n$ center problem (\ref{eq:ncent}), we consider also  systems with noncompact
configuration space $M$.  Then we have to impose certain conditions at infinity.
A simple way out is to assume that there is a compact geodesically convex domain in $M$
containing all singularities.

Let $D\subset M$ be a compact domain with smooth boundary $\partial D$.
When $M$ is a closed manifold, we   set $D=M$,
then $\partial D$ is empty.
Fix an energy level
\begin{equation}
\label{eq:D}
H=h>\max_{D}V.
\end{equation}
By Maupertuis' principle, trajectories $\gamma:[a,b]\to M$ with energy $h$ are (up to a reparametrization)
geodesics of the Jacobi metric
\begin{equation}
\label{eq:Jac}
\|\dot q\|_h=g_h(q,\dot q)=\max_p\{\langle p,\dot q\rangle: H(q,p)=h\}
=\sqrt{2(h-V(q))}\|\dot q\|,
\end{equation}
i.e.\ extremals of the Maupertuis--Jacobi action functional
\begin{equation}
\label{eq:J}
J(\gamma)=\int_a^b \|\dot \gamma(t)\|_h\, dt.
\end{equation}
Under condition (\ref{eq:D}),  $g_h^2$ is a positive definite Riemannian metric in $\hat D=D\setminus\Delta$.
To study the Hamiltonian flow on the level $\{H=h\}$ is the same as to study geodesics of the  Riemannian metric $g_h^2$.

The boundary $\partial D$
is called {\em geodesically convex} for energy $h$ if it is geodesically convex with respect
to the Jacobi metric.
Thus for every  trajectory $q(t)$ with energy $h$
such that $q(0)\in\partial D$ and $\dot q(0)\in T_{q(0)}(\partial D)$, there is $\eps>0$ such that
$q(t)\notin D\setminus \partial D$ for $-\eps<t<\eps$.
For example, a domain bounded by a non self-intersecting periodic trajectory of energy $h$ is geodesically convex.

Let $\nu$  be the inner unit normal vector to $\partial D$
with respect to the Riemannian metric $\|\cdot\|$,
and $\kappa$  the geodesic curvature
corresponding to the normal $\nu$. Thus if $\tau$ is the unit tangent vector and $s$ the arc length along $\partial D$, 
then
$$
\frac{D\tau}{ds}=\kappa \nu.
$$

By (\ref{eq:Newton}), the boundary is geodesically convex if for the motion with energy $h$ along the boundary,
the normal force $F_{norm}=-\langle \nabla V,\nu\rangle$ is smaller than the normal acceleration
$$
a_{norm}=\langle\frac{D\dot q}{dt},\nu\rangle= \kappa \|\dot q\|^2=2\kappa(h-V).
$$
Thus $\partial D$ is geodesically convex for energy $h$ iff
\begin{equation}
\label{eq:convex}
\langle \nabla V(q),\nu(q)\rangle+ 2\kappa(q)(h-V(q))\ge 0,\qquad q\in \partial D.
\end{equation}

We divide singularities into classes  depending on their strength. Let
\begin{equation}
\label{eq:Ak}
A_k=2-2k^{-1},\qquad k=1,2,3,\dots,\infty.
\end{equation}
Then
$$
A_1=0,\quad A_2=1,\quad A_3=4/3,\quad A_4=3/2,\;\dots\;, \;A_\infty=2.
$$
We will see that singularities $a_j$ with $\alpha_j=A_k$ are regularizable \cite{Knauf}:
$$
\Delta_{reg}=\{a_j:\alpha_j=A_k,\; k=2,3,\dots\}.
$$
Let
$$
\Delta_k=\{a_j: A_k\le \alpha_j<A_{k+1}\}.
$$
Hence
$$
\Delta_1=\Delta_{weak},\quad \Delta_\infty=\{a_j: \alpha_j\ge 2\}=\Delta_{jac}\cup\Delta_{strong}.
$$
Set $n_k=\# \Delta_k$ and
$$
A(\Delta)=\sum_{2\le k\le \infty} n_kA_k=n_2+\frac 43 n_3+\frac 32 n_4+\cdots+ 2 n_\infty.
$$
Next  we formulate the main result of the paper.

\begin{thm}
\label{thm:top}
Let $D\subset M$ be a compact domain containing all  singularities,
and let $h>\max_DV$.
Suppose that   the boundary $\partial D$ is geodesically
convex  for energy $h$.
If
\begin{equation}
\label{eq:chi}
A(\Delta) >2\chi(D),
\end{equation}
then:
\begin{itemize}
\item
There exist an infinite number of noncontractible
noncollision periodic orbits with energy $H=h$ in $\hat D=D\setminus\Delta$
and an infinite number of heteroclinic orbits joining these periodic orbits.
\item
The flow on the energy level  $\{H=h\}\cap T^*\hat D$
has a compact chaotic invariant set with positive topological entropy.
\end{itemize}
\end{thm}

When all singularities are Newtonian (then $A(\Delta)=n$)
Theorem \ref{thm:top} is an old result,  see e.g.\ \cite{Bol:ncent,Bol:sing,Knauf}.
Another simple case is when  all singularities  are strong  (then $A(\Delta)=2n$), see \cite{Gordon,Ambr}.

For singularities with $1\le \alpha_j< 2$, we have $A(\Delta)\ge n$,
so for closed $M$ Theorem \ref{thm:top} implies Theorem \ref{thm:simple}.
However, also then Theorem \ref{thm:top}  gives a  stronger statement.

When all singularities are weak with $0<\alpha_j<1$, we have $A(\Delta)=0$,
and so condition (\ref{eq:chi}) is $\chi(D)<0$, as if there are no singularities.
Thus weak singularities  are ignored in Theorem \ref{thm:top}.

Since $A(\Delta)\le \sum_{j=1}^n\alpha_j$, condition (\ref{eq:chi}) is slightly stronger than (\ref{eq:sum}).
The assertion of Theorem \ref{thm:top} was proved in \cite{Bol-Koz:Izv} under the additional assumptions
that $D$ is homeomorphic to a plane domain and $ \Delta=\Delta_{reg}$: all singularities are regularizable \cite{Knauf}.
Then $A(\Delta)=\sum \alpha_j$,
so condition (\ref{eq:chi}) coincides with (\ref{eq:sum}).

One can partly describe symbolic dynamics in the chaotic set in Theorem \ref{thm:top}.

\begin{thm}
\label{thm:cover}
There exists a surface $X$ with boundary, a $K$-sheet smooth covering $\phi:X\to D\setminus(\Delta_{jac}\cup \Delta_{strong})$
branched over the  set $\Delta_{newt}\cup\Delta_{mod}$ of Newtonian
and moderate singularities, and  a smooth complete\footnote{Thus the corresponding metric space is complete.
Equivalently, geodesics exist till they exit through $\partial X$.}
Riemannian metric   on  $X$  such that:
\begin{itemize}
\item 
Projections to $D$ of minimal geodesics  on the universal covering of the surface $X$ are trajectories with energy $H=h$ having no collisions with $\Delta$, except maybe with regularizable singularities $\Delta_{reg}$.
\item
The  Euler characteristics
$$
\chi(X)=K\big(\chi(D)-\frac 12 A(\Delta)\big)
$$
is negative when (\ref{eq:chi}) holds.
\end{itemize}
\end{thm}

Theorem \ref{thm:top} follows from Theorem \ref{thm:cover} and classical properties of geodesic flows on closed surfaces,
see e.g.\ \cite{KH}. Our surface $X$ has  a convex boundary, but main properties can be extended to this case,
see e.g.\ \cite{Bol:geo}.

\begin{rem} 
Using a formula in \cite{Katok} one can roughly estimate the topological entropy of the geodesic flow of the Jacobi metric:
$$
h_{top}\ge \sqrt{\frac{\pi(A(\Delta)-2\chi(D))}{\vol(D)}},
$$
where the volume is computed in the Jacobi metric.
\end{rem}

\begin{rem}
Theorem \ref{thm:top}  can be generalized to the case of systems with exact gyroscopic forces, 
when the gyroscopic 2-form $\omega$
in (\ref{eq:Omega}) is exact, i.e.\ it is a differential of a 1-form $\langle w(q),dq\rangle$.
By the change $p\to p=w(q)$ the symplectic form (\ref{eq:Omega})  can be replaced by the standard form $dp\wedge dq$, 
and the Hamiltonian by
$$
H(q,p)=\frac 12\|p-w(q)\|^2+V(q).
$$
Then condition (\ref{eq:h}) needs to be replaced by
$$
h>\max_{q\in D}\big(V(q)+\frac 12 \|w(q)\|^2\big).
$$
Under this condition the Jacobi metric
$$
g_h(q,\dot q)=\|\dot q\|_h+\langle w(q),\dot q\rangle
$$
is a positive definite Finsler metric on $D$.
Also the definition of geodesic convexity needs to be modified since the Jacobi metric is irreversible.
If the gyroscopic form is $\omega=u(q)\,\Omega$, where $\Omega$ is the area form on $M$,
we need to replace (\ref{eq:convex}) with
$$
\langle \nabla V(q),\nu(q)\rangle+ 2\kappa(q)(h-V(q))-|u(q)|\sqrt{2(h-V(q))}\ge 0,\qquad q\in \partial D.
$$
However for simplicity we consider only natural Hamiltonian systems. 
For systems with gyroscopic forces see \cite{Bol:sing}.
\end{rem}

\begin{rem}
\label{rem:libr}
As already mentioned,  Theorem \ref{thm:top} does not work for $h<\sup V$. Then only much weaker results can be proved.
Suppose that the domain of possible motion $D=D_h=\{q\in M:V(q)\le h\}$ is compact and $\nabla V\ne 0$ on $\partial D$.
The Jacobi metric vanishes on $\partial D$.  Theorem \ref{thm:cover} still works,
but now the metric in $X\setminus \partial X$ is of course incomplete. However, the results of \cite{Bol-Koz:libr}
imply that the number of minimizing geodesic arcs starting and ending on the boundary $\partial X$
is at least $\rank\, H_1(X,\partial X,\Z)$.
These geodesics are projected by the covering $\phi$ to  reversible periodic orbits (librations) 
with energy $h$  
which have no collisions with $\Delta$, except maybe regularizable singularities $\Delta_{reg}$.
If  there are no regularizable singularities, then there will be no collisions, 
so there exist at least $\rank\, H_1(X,\partial X,\Z)$
noncollision reversible periodic orbits with energy $h$.
In contrast to the case $h>\sup V$ we get only a finite number of periodic orbits and we have no hope
to get chaotic trajectories (at least by elementary topological methods of this paper).
\end{rem}

Theorems \ref{thm:top} and \ref{thm:cover} are proved in section \ref{sec:top}. 
In section \ref{sec:convex} we show that in the proof 
without loss of generality we may assume that there are no strong singularities.
In the main part of the paper we also assume that there are no Jacobi singularities.
In the presence of Jacobi singularities  the proof of the first part of Theorem \ref{thm:top} does  not change,
but the proof of the second part requires additional arguments.
The case of nonempty  $\Delta_{jac}$ is discussed in the last section.

\section{Examples}

1. $M=\T^2$ is a torus. Since $\chi(\T^2)=0$, by Theorem \ref{thm:top}, 
the existence of  $n\ge 1$  singularities with $\alpha_j\ge 1$
implies chaotic behavior
on  energy levels $h>\max V$. For Newtonian or strong singularities this is well known, 
see e.g.\ \cite{Bol:sing,Knauf,Knauf2}.
In \cite{Koz-Tre:sing} it is proved that  in the presence of singularities with $\frac12<\alpha_j<2$
on $\T^2$ there are no polynomial in momenta and integrable in coordinates  first integrals in the whole
phase space.\footnote{In \cite{Koz-Tre:sing} also systems with singularities on $\T^d$ were studied.}
But such nonintegrabilty in general does not imply chaotic behavior. We do not know if the existence
of a weak singularity with $\frac12<\alpha_j<1$ on $\T^2$ always implies positiveness of the topological entropy.

2. $M=S^2$ is a sphere. Then   $2\chi(S^2)=4$. If all singularities are strong, then $A(\Delta)=2n$,
 so Theorem \ref{thm:top} works for $n\ge 3$. A system on a sphere with 2 strong singularities
 may  be integrable.

Indeed, take a metric of revolution on a sphere,
place the singularities in the antipodal points and take the potential of revolution
depending only on the distance to these points. Then the angular
momentum about the axis of revolution will be a first integral.

If there are $n$ singularities with $\alpha_j\ge 3/2$, we have $A(\Delta)\ge \frac 32 n$,
so Theorem \ref{thm:top} works for $n\ge 3$. If there are $n$ singularities with $\alpha_j\ge 3/4$,
then  Theorem \ref{thm:top} works for $n\ge 4$.
If there are $n$ Newtonian singularities, then $A(\Delta)=n$, so
Theorem \ref{thm:top}   works for $n\ge 5$.
For $n=4$ Newtonian singularities the system may have a quadratic in the momentum first integral
on an energy level $H=h>\max V$,
see \cite{Bol:gyro}.

Indeed, take any Riemannian metric $\|\cdot\|$ on $S^2$ and any set $\Delta\subset S^2$ with $\#\Delta=4$.
The metric  defines on $S^2$ a conformal structure.
There exists a holomorphic differential $\omega$ of degree 2 having simple poles at each point in $\Delta$.
Then $|\omega|$ is a Riemannian metric on $S^2\setminus\Delta$ conformally equivalent to $\|\cdot\|$.
Define the kinetic energy by $T=\frac12\|\dot q\|^2$  and the potential energy so that $h-V=|\omega|/T$.
Then $V$ is a function
on $S^2\setminus\Delta$ and it has 4 Newtonian singularities.
The Jacobi metric $g_h^2$ is a constant multiple of $|\omega|$,
and so its geodesic flow  is integrable: it has a quadratic in momentum first integral $\Real\omega$.

If there are 3 Newtonian singularities
and the 4th a stronger one  with  $\alpha_4\ge \frac43$,  we get $A(\Delta)\ge 3+\frac43>4$,
so there is a chaotic invariant set for energies $h>\max V$.

3. The generalized $n$ center problem. Since $M=\R^2$ is noncompact, to apply Theorem \ref{thm:top} to the Hamiltonian  (\ref{eq:ncent})
 we need to  find a geodesically convex for energy $h>\sup V$ compact set  $D\subset \R^2$
containing all singularities.
We can try to take for $D$ a disk $B(0,R)=\{q:|q|\le R\}$
with sufficiently large radius $R$.
This  works under a convexity assumption (compare with (\ref{eq:convex})):
\begin{equation}
\label{eq:Dconv}
\langle\nabla V(q),q\rangle\le 2(h-V(q)),\qquad |q|= R.
\end{equation}
For example, (\ref{eq:Dconv}) holds  for large $R$  if $h>\sup V$ and
$$
\limsup_{|q|\to+\infty}\langle \nabla U(q),q\rangle\le 0.
$$
Then for large $R$ the disk $D= B(0,R)$ is geodesically convex for energy $h>\sup V$.

\begin{thm}
\label{thm:ncent}
If $A(\Delta)>2$,
the generalized $n$-center problem in $\R^2$ has a compact chaotic invariant set
and an infinite number of minimizing (on a suitable  branched cover)
periodic noncollision solutions on the energy level $H=h>\sup V$.
\end{thm}

We will show that Theorem \ref{thm:ncent} holds  without the convexity condition (\ref{eq:Dconv})
but then   the  choice of a geodesically convex compact domain $D\subset\R^2$  is less evident.
The proof is given in section \ref{sec:cone}.

Trajectories in the chaotic set  are projections of minimal geodesics of a smooth Riemannian metric
on a geodesically convex  compact $K$-sheet covering $X$ of the domain $D$ with
$$
\chi(X)=K\big(1-\frac 12 A(\Delta)\big)<0.
$$
If the degrees of all singularities belong to a single interval $A_k\le a_j<A_{k+1}$, 
i.e.\ $a_j\in \Delta_k$ for a single $k$, then the covering has degree $K=k$
and the surface $X$ is homeomorphic to the Riemannian surface of the function
$$
z\to \sqrt[k]{(z-a_1)\cdots (z-a_n)}.
$$
In general we may set
$$
K=\prod_{2\le k<\infty,\;n_k\ne 0}k,
$$
which of course is not optimal.

For example, if there is one Newtonian singularity   and another with $\alpha_2\ge 4/3$,
we get $A(\Delta)\ge 1+\frac 43>2$, so  Theorem  \ref{thm:ncent}
applies.   Theorem \ref{thm:ncent} also
gives chaotic behavior  if there are 3 or more   singularities at least of Newtonian strength.

In a recent paper \cite{Castel} it was proved that for the classical $n$-center problem with $U=0$ and 
$\alpha_1=\dots=\alpha_n\in (0,1)$,
minimizers of the Maupertuis length functional in  certain admissible
homotopy classes of closed curves in $\R^2\setminus \Delta_{mod}$ do not have collisions with singularities.
We prove a slightly more general proposition in section \ref{sec:cone}.
The  result of \cite{Castel} implies  chaotic behavior for the $n$ center problem  with $n\ge 4$  moderate singularities.
The assumption of  Theorem \ref{thm:ncent}
is considerably weaker. 

\begin{rem}
The case $U=0$ is  relatively simple because  then the Jacobi metric
has negative Gaussian curvature, see \cite{Knauf}.
Then Theorem \ref{thm:ncent} can be improved.
For example, it holds also when $h=\sup V=0$.
Indeed, then
$$
2(-V(q))- \langle\nabla V(q),q\rangle=\sum \frac{m_j}{|q-a_j|^{\alpha_j}}(2-\alpha_j+O(|q|^{-1})). 
$$
Hence if there exists at least one singularity with $\alpha_j<2$,  
then for large  $R$  the boundary of the disk
$B(0,R)$  is geodesically convex for energy $h=0$.
\end{rem}

\begin{rem}
Conditions of Theorem \ref{thm:ncent}   are purely topological:
no analytical properties except the presence of singularities are involved.
Of course under additional analytical assumptions much stronger results can be proved  also for the $n$-center problem in $\R^d$.
For example, for Newtonian singularities and large energy $h\to +\infty$ Aubry's method of antiintegrable limit can be used
to describe symbolic dynamics of chaotic orbits,
see e.g.\ \cite{Knauf,Knauf2,Bol-Mac:ncent,ANTI,Bol:MIAN}.
But this requires $n\ge 4$ centers not lying on a single line.
Theorem \ref{thm:top} works for $n\ge 3$ and all energies $h>\sup V$.
\end{rem}

\begin{rem}
We conjecture that  Theorem \ref{thm:ncent}  holds also for the $n$ center problem  in $\R^3$.
If all singularities are Newtonian,  then this was proved in
\cite{Bol-Neg:reg} by using the KS regularization and  the results of Gromov and Paternain, see \cite{Paternain}.
Elementary methods of the present paper will not work since  they use that $\R^2\setminus\Delta$
is not simply connected.
\end{rem}

In the rest of the paper we prove Theorems \ref{thm:simple} and \ref{thm:top}, 
and also remove the convexity assumption in Theorem \ref{thm:ncent}.

\section{Convexity properties}

\label{sec:convex}

 First we  prove   simple convexity properties of small  neighborhoods of singularities.
 Let $a_j\in\Delta$ be   a singularity of order $\alpha_j$
 and let $B_j=B(a_j,\eps)$ be the  small disk (\ref{eq:ball}).

 \begin{lem}
 \label{lem:convex}
 For $0<\alpha_j<2$ and small $\eps>0$
 the boundary circle $S_j=\partial B_j$ is geodesically convex in the Jacobi metric.
 For   $\alpha_j>2$ and small $\eps>0$
 the boundary circle $S_j$ is geodesically concave in the Jacobi metric, i.e.\
 the complement $\overline{D\setminus B_j}$ is geodesically convex.
 \end{lem}

 This property was   known already to the founders of celestial mechanics (it follows from the Lagrange--Jacobi identity).
 However, for completeness we give a proof.  For any $x\in S_j$, let $\nu(x)$ be the inner unit normal vector
 and $\kappa(x)$ the curvature (with respect to the metric $\|\cdot\|$). Then by (\ref{eq:V}),
 \begin{eqnarray*}
 \kappa(x)&=&\eps^{-1}+o(\eps^{-1}),\\
 V(x)&=&-m_j\eps^{-\alpha_j} +o(\eps^{-\alpha_j}),\\
 \nabla V(x)&=&-\alpha_j m_j\eps^{-\alpha_j-1}\nu(x)+  o(\eps^{-\alpha_j-1}).
 \end{eqnarray*}
 Hence
 \begin{eqnarray*}
 2(h-V(x))\kappa(x)&=&2m_j\eps^{-\alpha_j-1} +o(\eps^{-\alpha_j-1}),\\
 -\langle \nu,\nabla V(x)\rangle&=&\alpha_j m_j\eps^{-\alpha_j-1}+  o(\eps^{-\alpha_j-1}).
 \end{eqnarray*}
 Now the conclusion follows from (\ref{eq:convex}).
 \qed

 \medskip

 For  critical Jacobi singularities with $\alpha_j=2$ Lemma \ref{lem:convex} does not work.
 Such singularities have to be treated separately.

 \begin{cor}\label{cor:strong}
 In the proof of Theorem \ref{thm:top}
 without loss of generality we may assume that
 there are no strong singularities.
 \end{cor}

 By Lemma \ref{lem:convex}, for a strong singularity   $a_j$  with $\alpha_j>2$, the disk   $B_j=B(a_j,\eps)$
 has concave boundary $\partial B_j$. Hence the complement
 $$
 D'=\overline{D\setminus\cup_{\alpha_j>2}B_j}
 $$
 is a geodesically convex compact domain containing  the set $$
 \Delta'=\Delta\setminus\Delta_{strong}=\{a_j\in\Delta:\alpha_j\le 2\}
 $$
 of not strong singularities.
 We have
\begin{eqnarray*}
 A(\Delta')&=& \sum_{2\le k<\infty} n_kA_k=A(\Delta)-2n_\infty,\\
  \chi(D')&=&\chi(D)-n_\infty.
\end{eqnarray*}
 Hence
 $$
 A(\Delta)-2\chi(D)=A(\Delta')-2\chi(D').
$$
Replacing $D$ with $D'$ and $\Delta$ with $\Delta'$ in Theorem \ref{thm:top}
we get a  system without strong singularities satisfying the conditions of Theorem \ref{thm:top}.
\qed

\section{Jacobi's metric space}

\label{sec:cone}

Trajectories with energy $h$ are
geodesics of the Jacobi metric (\ref{eq:Jac}),
i.e.\ extremals of the Maupertuis--Jacobi action functional (\ref{eq:J}).

\begin{rem}
In many recent papers (see e.g.\ \cite{Castel,Terra1,Terra2}) a different Maupertuis functional is used:
$$
I(\gamma)=\left(\int_a^b\|\dot\gamma(t)\|^2\,dt\right)\left(\int_a^b(h-V(\gamma(t)))\,dt\right).
$$
Then $I(\gamma)\ge J(\gamma)^2/2$ and the equality holds iff the energy is constant along $\gamma$.
The functional $I$ has an advantage of being differentiable on an appropriate Sobolev space.
Using this functional makes sense when looking for minimax geodesics. However, if, as in the present paper,
only minimal geodesics are studied, then the classical Maupertuis functional $J$ in Jacobi's form
is more convenient since we can replace Analysis by simple metric geometry.
\end{rem}

Since the length $J(\gamma)$ of a curve is independent of a parametrization, we identify curves which
differ by an orientation preserving reparametrization and parametrize all curves by the segment $[0,1]$.
The Jacobi metric $g_h$ makes $\hat D=D\setminus\Delta$  a metric space with the distance
 $$
 \rho(x,y)=\inf\{J(\gamma):\gamma\in C^1([0,1],\hat D),\;\gamma(0)=x,\; \gamma(1)=y\}.
 $$
 If $a_j$ is a strong singularity with $\alpha_j>2$, then as $y\to a_j$, we have
 $$
 \rho(x,y)\sim d(y,a_j)^{1-\alpha_j/2}\to +\infty.
 $$
 For Jacobi singularities with $\alpha_j=2$,
 $$
 \rho(x,y)\sim -\ln d(y,a_j)\to +\infty.
 $$
 For singularities with $\alpha_j<2$,  there exists the limit
 $$
 \rho(x,a_j)=\lim_{y\to a_j}\rho(x,y)<+\infty.
 $$

Thus the metric is complete at Jacobi and strong singularities.
This was discovered by Poincar\'e and first used by Gordon \cite{Gordon}.
Hence the completion of the metric space $(\hat D,\rho)$
is the complete metric space $(D\setminus(\Delta_{jac}\cup\Delta_{strong}),\rho)$
with the standard topology.

\medskip

From now on we assume that there are no Jacobi singularities with  $\alpha_j=2$.
The   critical case when some $\alpha_j=2$  is treated separately in the last section.
By Corollary \ref{cor:strong}, we can assume that there are no strong singularities, so $0<\alpha_j<2$ for all $j$:
we have only weak, moderate and Newtonian singularities.
Then the completion of the metric space $(\hat D,\rho)$ is the compact metric space $(D,\rho)$.

Now we can define the length $J(\gamma)\in [0,+\infty]$ of any  curve
$\gamma\in C^0([0,1], D)$ in a standard way:
$$
J(\gamma)=\sup\Big\{\sum_{i=1}^k \rho(\gamma(t_{i-1}),\gamma(t_i))\;\Big|\;
0=t_0<t_1<\cdots<t_k=1\Big\}.
$$
It is natural to consider $J$ on the set
of rectifiable curves $\{\gamma\in C^0([0,1],D)\mid  J(\gamma)<\infty\}$
such that, if $\gamma$ is not a point curve, i.e.\ $J(\gamma)\ne 0$, then $\gamma(t)\ne\const$ on any interval.
Since curves differing by an order preserving reparametrization   are identified,
every element in  the corresponding quotient space $\tilde C([0,1],D)$ is uniquely represented
by a curve $\gamma:[0,1]\to D$ parameterized proportionally to the arc length:
$J(\gamma|_{[0,s]})= sJ(\gamma)$.
Then $\tilde C([0,1],D)$ is embedded in $C^0([0,1],D)$ and carries the $C^0$ topology.

It is well known \cite{Metric} that the length functional $J$ is lower semicontinuous
and for any $c>0$, the set
$\{\gamma\in \tilde C([0,1],D) : J(\gamma)\le c\}$ is compact.
Then $(D,\rho)$ is a compact length space \cite{Metric}:
for any points $x,y\in D$ there is a minimizing curve $\gamma$
joining them such that $J(\gamma)=\rho(x,y)$. By convexity,
if $x,y\in D\setminus\partial D$, the minimizer will not touch $\partial D$.

Let $\Gamma\subset C^0([0,1],\hat D)$ be a homotopy class of curves joining two points $x,y\in\hat D$.
The compactness property implies

\begin{lem}
\label{lem:hom}
The functional $J$ has a minimum on the closure  $\bar\Gamma$ of $\Gamma$ in $C^0([0,1],D)$.
Any minimizer $\gamma$ is a  trajectory with energy $h$ or  a chain of
trajectories $\gamma=\gamma_1\cdots\gamma_k$,
where the curves $\gamma_i$ join pairs of   singular points in $\Delta$.
\end{lem}

The same   holds if $\Gamma$ is a nontrivial homotopy class of closed curves in $\hat D$, i.e.\
a path connected component in $C^0(S^1,\hat D)$, where $S^1=[0,1]/\{0,1\}$.
We call a homotopy class of closed curves in $\hat D$
nontrivial if it does not  contain contractible loops or loops $\gamma:S^1\to B_j$   in small neighborhoods $B_j$
of singularities $a_j$. For a trivial class $\Gamma$, the minimum of $J$ on $\bar\Gamma$ is attained 
on a point curve $\gamma=a_j$.

In order to get noncollision trajectories with energy $h$ we need to show that the minimizer
$\gamma$ does not pass through the singular set $\Delta$.

Let $a_j\in\Delta$, $0<\alpha_j<2$. Take sufficiently small $\eps>0$ and let
$B_j=B(a_j,\eps)$ be the closed disk in Lemma \ref{lem:convex}
and $S_j$ the corresponding circle.
Since the disk $B_j$ is geodesically convex in the Jacobi metric,
a minimizer  joining a pair of points  $x,y\in B_j$ stays in $B_j$.
In particular,
$$
\inf\{J(\gamma)\mid \gamma\in  C^0([0,1],B_j),\;\gamma(0)=x,\;\gamma(1)=y\}=\rho(x,y).
$$

The next lemma means  that the metric $\rho$ has  a cone singularity at $a_j$
with total angle less than $2\pi$.

\begin{lem}
\label{lem:2pi}
Let $0<\alpha_j<2$. Then there exists $\lambda\in (0,1)$ such that for sufficiently small $\eps>0$ and
any points $x,y\in S_j$ we have
\begin{equation}
\label{eq:rho}
\rho (x,y)<\lambda(\rho(x,a_j)+\rho(a_j,y)).
\end{equation}
\end{lem}

We will prove Lemma \ref{lem:2pi} in the next section.   By (\ref{eq:rho}),
$$
\rho (x,y)< \rho(x,a_j)+\rho(a_j,y)-2r,\qquad r= (1-\lambda)\rho(S_j,a_j).
$$
If a  curve $\gamma$ joining $x,y\in S_j$ enters  the ball
$$
\calB_j=\calB(a_j,\delta)=\{x:\rho(x,a_j)\le \delta\}\subset B_j,\qquad 0<\delta<r,
$$
then
$$
J(\gamma)\ge \rho(x,a_j)+\rho(a_j,y)-2\delta> \rho(x,y)+\mu,\qquad \mu=2r-2\delta.
$$
Thus if
$$
J(\gamma)\le   \rho(x,y)+\mu,
$$
then $\gamma$ does not enter the  ball $\calB_j$.

\begin{cor}
\label{cor:weak}
Let $\pi:\tilde D\to D$ be the universal covering, $\tilde\Delta=\pi^{-1}(\Delta)$, $\tilde\rho$
the corresponding  distance on $\tilde D$, and $\tilde J$ the length functional.
If $\gamma$ is a minimizer of $\tilde J$
joining a pair of points in $\tilde D\setminus \tilde \Delta$,
then $\gamma$ does not pass though   $\tilde\Delta$.
\end{cor}

This implies the following weak version of Theorem \ref{thm:top}.

\begin{prop}\label{prop:weak}
Suppose that all singularities satisfy $0<\alpha_j<2$.
Every nontrivial homotopy class of closed curves
in $D$ contains a minimal noncollision closed geodesic.
If $\chi(D)<0$ then there exist an infinite number of minimal heteroclinic
geodesics. The geodesic flow of the Jacobi metric has a compact chaotic invariant set
of minimal (on the universal covering of $D$) noncollision geodesics
with positive topological entropy. 
\end{prop}

When there are no singularities ($\Delta=\emptyset$) and no boundary ($D=M$ is a closed surface of genus $\ge 2$),
this is an old result essentially known
to Morse \cite{Morse} and Hedlund \cite{Hedlund} (of course except the definition of the topological entropy).
They proved the existence of many minimal heteroclinic geodesics
joining minimal closed geodesics in homotopy classes. In fact the set of
minimizing geodesics on the universal covering $\tilde M$
after projecting to $M$ gives a compact invariant set  for the geodesic flow on $M$
with  positive topological entropy \cite{KH}.
This is true also for surfaces $D$ with geodesically convex boundary,  see e.g.\ \cite{Bol:geo}.

Indeed, suppose  that $\partial D$ consists of closed curves $C_1,\dots,C_n$.
Take a minimal closed geodesic $\gamma_i$ in the homotopy class  of the  curve $C_i$ 
(for example by applying the curvature flow to $C_i$).
Then $\gamma_i\cap \gamma_j=\emptyset$  for $i\ne j$, so $S=\cup \gamma_i$ 
bounds a domain $U\subset D$ homeomorphic to $D$ and with zero curvature of the boundary. 
Glue two copies of $U$
along $S$ and obtain a closed manifold $N$ without boundary. Since the boundary
curves $\gamma_i$  
have zero curvature,
one can show that the Riemannian metric on $U$ defines a $C^2$ Riemannian metric $\bar g$ on $N$
invariant under the involution $\sigma:N\to N$ interchanging the copies of $U$.
Now the results of Morse and Hedlund on minimal heteroclinic geodesics can be applied \cite{Bol:geo}.

If there are only weak singularities with $0<\alpha_j<1$, then Proposition \ref{prop:weak}
coincides with Theorem \ref{thm:top}.
In fact we will  deduce Theorem \ref{thm:top} from
Proposition \ref{prop:weak} by using successive Levi-Civita type regularizations.

\medskip

\noindent{\it Proof of Proposition \ref{prop:weak}.}
Let us modify the Jacobi metric $g_h$
in every ball $\calB_j$  replacing it by a smooth Riemannian metric $g'\ge g_h$ on $D$
such that $g'=g_h$
on $D\setminus \cup \calB_j$. Then every nontrivial homotopy class of closed curves in $D$
contains a minimal geodesic, there is an infinite number of heteroclinic geodesics
joining minimal geodesics, and there is a compact chaotic invariant set
of minimizing (in the universal covering $\tilde D\to D$) geodesics of the metric $g'$.
By Lemma \ref{lem:2pi} such geodesics can not pass through the balls $\calB_j$,
so they are geodesics of the  Jacobi metric $g_h$.
\qed

\medskip

For moderate singularities with $1<\alpha_j<2$  a stronger version of Lemma \ref{lem:2pi} is true:
the total angle at the vertex $a_j$ of the cone is less than $\pi$.
Let $\phi:B'_j\to B_j$ be a 2-sheet covering branched over $a_j$ (same as in the Levi-Civita regularization,
see the next section). Thus $a'_j=\phi^{-1}(a_j)$ is a single point and for $x\ne a_j$, $\phi^{-1}(x)$ consists of two points.
Lift the metric $\rho$ to a metric $\rho'$ on $B'_j$.

\begin{lem}
\label{lem:pi}
Let $1<\alpha_j<2$.
Then there is $\lambda\in (0,1)$ such that for sufficiently small $\eps>0$
and any $x,y\in S'_j=\partial B'_j$, we have
$$
\rho'(x,y) < \lambda (\rho'(x,a'_j)+\rho'(a_j',y)).
$$
\end{lem}

As in Lemma \ref{lem:2pi}, there exist $\delta,\mu>0$ such
any  curve $\gamma$ joining $x,y\in S'_j$
such that
$$
J'(\gamma)\le   \rho'(x,y)+\mu
$$
does not enter the  ball $\calB'_j=\{x:\rho'(x,a'_j)\le \delta\}\subset B'_j$.

\begin{cor}
\label{cor:mod}
Let $1<\alpha_j<2$.
Then a minimizer joining any points $x,y\in S'_j$
does not pass though  $\calB'_j$.
\end{cor}

Note that for Newtonian singularities this is not true: if $x\ne y$ are points on $S'_j$ with $\phi(x)=\phi(y)$,
then $\rho'(x,y)=\rho'(x,a'_j)+\rho'(a'_j,y)$.

Let  $\Delta_{mod}$ be
the set of moderate singularities with $1<\alpha_j<2$, and $\Delta_{newt}$ the set of Newtonian singularities.
Suppose  that $\Sigma=\Delta_{mod}\cup\Delta_{newt}\ne\emptyset$ (or else Theorem \ref{thm:top}
is reduced to Proposition \ref{prop:weak}).
The following proposition is an elementary property of Riemannian surfaces.\
It was used in \cite{Bol:sing} to deal with Newtonian singularities, see also \cite{Knauf} for $D=\T^2$.
We recall the proof  in section \ref{sec:simple} while proving Theorem \ref{thm:simple}.

\begin{prop}\label{prop:reg}
There exists a compact surface $D'$ with boundary and a  $K$-sheet ($K=2$ or $K=4$) smooth covering
$\phi:D'\to D$  
branched over $\Sigma$, so that each point in $\Sigma'=\phi^{-1}(\Sigma)$
has multiplicity 2. Thus $\#\phi^{-1}(a_j)=K/2$ for $a_j\in\Sigma$ and
$\#\phi^{-1}(x)=K$ for $x\in D\setminus \Sigma$.
\end{prop}

Lemma \ref{lem:pi} implies:

\begin{cor}
Let $\pi:\tilde D\to D'$ be the universal covering of $D'$, $\Phi=\phi\circ \pi:\tilde D\to D$, 
and  $\tilde\Delta=\Phi^{-1}(\Delta)$.
Then a minimizer joining a pair of points in $\tilde D\setminus\tilde\Delta$
does not pass though $\tilde\Delta$, except maybe regularizable Newtonian singularities $\Phi^{-1}(\Delta_{newt})$.
\end{cor}

Lemma \ref{lem:pi} will be proved  in the next section by
using the generalized Levi-Civita regularization.
It admits the following reformulation.

For given $x,y\in S_j$ there exist 2 simple homotopy classes $\Gamma_\pm(x,y)$ of curves
in $C^0([0,1], B_j\setminus\{a_j\})$ joining  $x,y$ while passing on different sides of $a_j$. Let
 $$
 \rho_\pm(x,y)=\inf\{J(\gamma):\gamma\in \Gamma_\pm(x,y)\}.
 $$
 Then of course   $\rho(x,y)=\min (\rho_+(x,y),\rho_-(x,y))$.

\begin{lem}\label{lem:2pi2}
Let $1< \alpha_j<2$. Then there is $\lambda\in (0,1)$ such that for sufficiently small $\eps>0$ and any $x,y\in S_j$,
we have
$$
\max (\rho_+(x,y),\rho_-(x,y))< \lambda (\rho(x,a_j)+\rho(a_j,y)).
$$
\end{lem}

As for Lemma \ref{lem:2pi} this implies that there exist $\delta,\mu >0$ such that any  curve $\gamma\in\Gamma_\pm(x,y)$
such that $J(\gamma)\le   \rho_\pm(x,y)+\mu$ does not enter the  ball $\calB_j$.

As an application, we give a simple proof of the main result of \cite{Castel}.
Let $\Sigma=\Delta_{newt}\cup\Delta_{mod}$.
Following \cite{Castel}, we call a homotopy class $\Gamma$ of closed curves in $C^0(S^1,D\setminus \Sigma)$
admissible if its  representative with a minimal number of self intersections
has no simple subloops which bound a disk containing a single singularity $a_j\in\Sigma$.

\begin{cor}
Let $\Gamma$ be an admissible   homotopy class  of closed curves in $C^0(S^1,D\setminus \Sigma)$.
Then any minimizer $\gamma$ of $J$ in the closure of $\Gamma$ in $C^0(S^1,D)$
has no collisions with moderate singularities $\Delta_{mod}$. If $\gamma$ has a collision
with $\Delta_{newt}$ at $t=\tau$, then it is reversible:
$\gamma(\tau+t)=\gamma(\tau-t)$.
\end{cor}

Indeed, suppose that $\gamma(\tau)=a_j\in\Delta_{mod}$.
There  exists a minimizing sequence $\gamma_k\to\gamma$ with a minimal number of self intersections, see \cite{Top}.
Then $\gamma_k$ enters the ball $\calB_j$ for large $k$,
so there is a segment $C$ of $\gamma_k$ joining the points $x,y\in S_j$
 in $B_j$ and crossing $\calB_j$. If the class $\Gamma$ is admissible, then $C$
 does not contain a subloop going around $a_j$, so $C\in\Gamma_+(x,y)$ or $C\in\Gamma_-(x,y)$.
 Then by Lemma \ref{lem:2pi2}, $J(C)>\rho_\pm(x,y)+\mu$, so $J(\gamma)>\inf_{\Gamma}J+\mu$, 
 a contradiction for large $k$.
\qed

\medskip

Now we use  Lemma \ref{lem:2pi2} to prove Theorem \ref{thm:ncent} without the convexity assumption.

\medskip

\noindent{\it Proof of   Theorem \ref{thm:ncent}.}
Take a simple closed curve in $\R^2$  encircling all singularities and minimize  the length functional $J$
on the homotopy class $\Gamma$
of this curve in $\R^2\setminus (\Delta_{newt}\cup \Delta_{med})$. Since
$$
g_h(q,\dot q)\ge c |\dot q|,\qquad 0<c<\sqrt{2(h-\sup V)},
$$
by Lemma \ref{lem:hom}, the minimum will be achieved on a closed curve $\gamma$ in the closure $\bar\Gamma$
of the class $\Gamma$ in $C^0(S^1,\R^2)$.
The curve $\gamma$ has no self intersections.
Hence by Lemma \ref{lem:2pi2}, $\gamma$ can not pass through the singular set $\Delta_{mod}$. If $\gamma$
passes through $\Delta_{newt}$ then $\gamma$ will be a segment joining 2 Newtonian singularities
passed twice in the opposite directions. Then there are no other non weak
singularities. This is impossible since $A(\Delta)>2$. 

By Lemma \ref{lem:pi}, $\gamma$ can not pass also though weak singularities. Hence $\gamma$ bounds
a compact domain $D\subset\R^2$ with geodesically convex boundary containing all Newtonian and moderate singularities.
Now Theorem \ref{thm:ncent} follows from Theorem \ref{thm:top}.
\qed

\section{Levi-Civita regularization}

\label{sec:Levi}

The Levi-Civita regularization is the main tool of this paper.

In  a neighborhood of any point of $M$ there exist local
coordinates $q_1,q_2$ such that the  Riemannian metric defining the kinetic energy has a conformal form:
\begin{equation}\
\label{eq:conf}
\|\dot q\|^2=g(q_1,q_2)(\dot q_1^2+\dot q_2^2)=g(z)|\dot z|^2,
\end{equation}
where  $z=q_1+iq_2\in\mC$.
Passing to an orienting 2-sheet covering we may assume that $M$ is oriented.
Then conformal  coordinates endow $M$ with a structure
of a Riemannian surface.

Let us choose conformal coordinates in a neighborhood of a singular point $a_j\in\Delta$ so that it corresponds to   $z=0$.
After a conformal change of variables without loss of generality we may assume that in
the  metric  (\ref{eq:conf}), $g(z)=1+O(|z|^2)$.
Then for the distance in the metric $\|\cdot\|$ we have
$$
d(q,a_j)=|z|(1+O(|z|^2)).
$$
By (\ref{eq:V}), the potential energy in a neighborhood of the point
$a_j$ takes the form
\begin{equation}
\label{eq:V(z)}
V=-\frac{f(z)}{|z|^\alpha}(1+O(|z|^2))+U(z),\qquad \alpha=\alpha_j,\quad f(0)=m_j>0.
\end{equation}

Following the idea of the Levi-Civita regularization \cite{L-C}, let us make a change of variables
\begin{equation}
\label{eq:Levi}
z=\phi_\beta(w)=w^\beta,\qquad w\in B(0,\eps)=\{w\in\mC:|w|\le\eps\},
\end{equation}
where the real number $\beta>1$ have to be determined. For integer $\beta=k$
the change $\phi_k:B(0,\eps)\to B(0,\eps^k)$ is a smooth $k$-sheet covering
branched at $w=0$. For noninteger  $\beta$ the map  $\phi=\phi_\beta$ is correctly defined  
by the formula
$$
\phi(w)=r^\beta e^{i\beta\theta}
$$ 
in  the domain
$$
\{w=r e^{i\theta}:\; 0\le r\le  \eps,\; 0\le \theta<2\pi\}.
$$

Instead of making a change of variables (\ref{eq:Levi}) in the Hamiltonian, we will work with the Jacobi metric
$$
g_h^2=2(h-V(z))g(z)|\dot z|^2
$$
corresponding to the energy level $H=h>\sup V$.
By (\ref{eq:V(z)}), in the conformal coordinates near the singular point $a_j$, the Jacobi metric is
\begin{equation}
\label{eq:gsing}
g_h^2=2g(z)\left(\frac{f(z)}{|z|^{\alpha}}(1+O(|z|^2))+h-U(z)\right)|\dot z|^2.
\end{equation}
After the  transformation $\phi$, the   metric $\tilde g=\phi^* g_h$  takes the form
\begin{equation}
\label{eq:tildeg}
\tilde g^2=
2\beta^2g(w^\beta)\left(\frac{f(w^\beta)}{|w|^{\tilde\alpha}}(1+O(|w|^{2\beta}))+(h-U(w^\beta))|w|^{2(\beta-1)}\right)|\dot w|^2,
\end{equation}
where
\begin{equation}
\label{eq:tilde_alpha}
\tilde\alpha=\alpha\beta-2(\beta-1)=2-\beta(2-\alpha).
\end{equation}

Let us choose $\beta$ so that  $\tilde\alpha=0$. Then
\begin{equation}
\label{eq:beta}
\alpha=2-\frac{2}{\beta},\quad \beta=\frac{2}{2-\alpha}.
\end{equation}
If $0< \alpha<2$, then $1< \beta<\infty$ (and vice versa).
The Jacobi metric  takes the form
$$
\tilde g^2= 2\beta^2g(w^\beta)(f(w^\beta)(1+O(|w|^{2\beta}))+(h-U(w^\beta))|w|^{2(\beta-1)})|\dot w|^2.
$$

If $\beta=k$ is an integer, then $\alpha=A_k$ (see (\ref{eq:Ak})), and
the   metric $\tilde g^2$ is  smooth  and positive definite on the disk $B(0,\eps)$.
Then $\phi=\phi_k$ is the generalized Levi-Civita regularization of order $k$ introduced in \cite{Knauf}.
 For Newtonian singularities with $\alpha=1$, $\beta=2$ we have the classical Levi-Civita
 regularization $z=w^2$.

Equations (\ref{eq:tildeg})--(\ref{eq:tilde_alpha}) imply:

\begin{lem}
\label{lem:LC}
Let   $A_k< \alpha_j<A_{k+1}$. Then the generalized Levi-Civita transformation of order $k$ centered at $a_j$
transforms the singularity $a_j$
to a weak singularity of order $\tilde\alpha_j=2-k(2-\alpha_j)$. Singularities with $\alpha_j=A_k$ will disappear.
Jacobi singularities are transformed  to  Jacobi singularities with $\tilde\alpha_j=\alpha_j=2$.
\end{lem}

 \begin{rem}
For Jacobi singularities,  the substitution $z=e^w$ is more natural. However, we will not use it in this paper.
\end{rem}

We are interested in the case of noninteger $\beta$ when regularization is not possible.
Let $W$ be the cone obtained from the sector
$$
\{w=r e^{i\theta}:\; 0\le r\le  \eps,\; 0\le \theta\le 2\pi/\beta\}
$$
 by identifying $ \theta=0$ with $\theta=2\pi/\beta$.
 The total angle at the vertex $w=0$ of the cone is $2\pi/\beta<2\pi$.
Then
$$
\phi:W\to B(0, \eps^\beta)
$$
is a  homeomorphism, and a diffeomorphism away from   the vertex.
The metric is smooth in the cone $W$ except  at the vertex $w=0$,
and  nearly Euclidean for $\beta>1$:
\begin{equation}
\label{eq:c}
\tilde g=(c+O(|w|^{2(\beta-1)}))|\dot w|,\qquad c=\beta\sqrt{2f(0)g(0)}.
\end{equation}
Let $L(\gamma)$ be the length of a curve $\gamma$ in $W$ in the Euclidean metric $|\dot w|$, 
and $J(\gamma$ in the Jacobi metric $\tilde g$. Then
\begin{equation}
\label{eq:L-J}
|J(\gamma)-c L(\gamma)|\le C\eps^{2(\beta-1)} L(\gamma).
\end{equation}

\medskip\noindent {\it
Proof of Lemma \ref{lem:2pi}}.
Since $\beta>1$, the angle at the vertex $w=0$ of the cone $W$ is  $2\pi/\beta<2\pi$.  Hence  for any  points
 $x,y\in \partial W$ the Euclidean distance in $W$ is not more than the length $2\eps\sin(\pi/2\beta)$ of a chord
 of the circle $|w|=\eps$ with angle $\theta=\pi/\beta<\pi$. Hence
$$
d(x,y)\le 2\eps\sin(\pi/2\beta)<2\eps= d(x,0)+d(0,y).
$$
Take $\lambda\in (\sin(\pi/2\beta),1)$. Then
\begin{equation}
\label{eq:dist}
d(x,y)\le  \lambda(d(x,0)+d(0,y)).
\end{equation}
By (\ref{eq:L-J}), for small $\eps>0$ a similar estimate holds for the metric $\rho$.
\qed

\medskip\noindent {\it
Proof of Lemma \ref{lem:pi}}. If $1<\alpha<2$, then $\beta>2$. Let $W'$ be the cone obtained from
$$
\{w=r e^{i\theta}: 0\le r\le \eps,\; 0\le \theta\le 4\pi/\beta\}
$$
by identifying
$\theta=0$ with $\theta=4\pi/\beta<2\pi$. Then we have a double covering $\phi:W'\to B(0,\eps^\beta)$ which
 is smooth except at the vertex. The angle at the vertex of $W'$ is less than $2\pi$.
Thus exactly the same  estimate for the Euclidean distance between points $x,y\in W'$ gives
(\ref{eq:dist}) with $\lambda\in (\sin(\pi/\beta),1)$.
\qed

\medskip

For Newtonian singularities, the angle at the vertex of the cone $W'$ will be $2\pi$, so the proof fails.
Then the singularity completely disappears after the Levi-Civita regularization, so
minimizing trajectories in $K'$ may pass through the singularity.

\section{Proof of Theorem \ref{thm:simple}}

\label{sec:simple}

 We follow \cite{Bol:sing}, where only Newtonian singularities were studied.

\begin{lem}
\label{lem:even}
Let $\Lambda=\{a_1,\dots,a_k\}\subset D$ be a finite set.
Suppose $k=\#\Lambda$ is even or   $D$ is homeomorphic to a domain in the plane.
There exists a  smooth double covering $\phi:D'\to D$, branched over $\Lambda$,
such that near each   point in $\Lambda'=\phi^{-1}(\Lambda)$ the map
 $\phi$ is the classical Levi--Civita regularization.
\end{lem}

\noindent{\it Proof.}
First suppose that $D$ is conformally equivalent to a  domain  in the complex plane.
Let $Z$ be the hyperelliptic Riemannian surface
\begin{equation}
\label{eq:Z}
Z=\{(z,w)\in \mC^2:   w^2=\prod_{j=1}^k(z-a_j)\}.
\end{equation}
Then $Z$ is a smooth  surface and the projection $\pi:Z\to \mC$, $\pi(z,w)=z$, is
a branched covering which locally near $a'_j=(a_j,0)=\pi^{-1}(a_j)$ is  the classical Levi-Civita transformation
$z=a_j+ c_j w^2+\cdots$. We set $D'=\pi^{-1}(D)\subset Z$ and $\phi=\pi|_{D'}$.

If $D$ is not homeomorphic to a plane domain,  take a domain $U\subset D\setminus \partial D$
conformally equivalent to a disk in $\mC$
 such that $\Lambda\subset U$ and define a surface $U'=\pi^{-1}(U)\subset Z$ and a branched double covering $\phi:U'\to U$
 as  above with $D$ replaced by $U$.
Next we extend $\phi:U'\to U$ to a double  covering $\phi:D'\to D$ as follows.

The boundary $\partial U$ is a closed curve. If $k=\#\Sigma$  is even,
$\partial U'=\phi^{-1}(\partial U)$ consists of 2 components -- closed curves $S_1,S_2$,
and $\phi:S_i\to \partial U$ is a diffeomorphism (for odd $k$, the boundary $\partial U'$
consists of a single closed curve and $\phi:\partial U'\to \partial U$ is a double covering).
We attach  to $U'$ two copies of $D\setminus U$ along the boundary circles $S_1,S_2$  and obtain
a  smooth surface $D'$ and a smooth double covering $\phi:D'\to D$ branched over $\Lambda$.
\qed

\medskip

We already removed all strong force singularities, and assumed that there are no Jacobi singularities.
Let $\Sigma=\Delta_{newt}\cup\Delta_{mod}$.
Let us prove Theorem \ref{thm:simple} for the case when $k=\#\Sigma$ is even or $D$ is a plane domain.
Set $\Lambda=\Sigma$ in Lemma \ref{lem:even}.
We lift the Jacobi metric $g_h$ on $D$ to a Riemannian metric  on $D'$ with the singular set $\Delta'=\phi^{-1}(\Delta)$.
For this metric, Newtonian singularities  will
disappear by the classical result of Levi-Civita.
Weak singularities will remain weak,
but their number will double. By Lemma \ref{lem:LC},
moderate singularities with $1<\alpha_j<3/2$ will become weak
of order $\alpha'_j=2\alpha_j-2\in (0,1)$,
while singularities with $3/2< \alpha_j< 2$ will become  moderate. Singularities with $\alpha_j=3/2$
will become Newtonian with $\alpha'_j=1$.

By the Riemann--Hurwitz formula,
$$
\chi(D')=2\chi(D)-k<0.
$$
Now the assertion of Theorem \ref{thm:simple} follows from Proposition \ref{prop:weak}.

\medskip

Lemma \ref{lem:even} implies Proposition \ref{prop:reg} for even $k=\#\Sigma$.
If  $k$ is odd, we first take a  covering of $D$ making it even.
If $D$ is not simply connected,  there exists a double covering $\psi:\tilde D\to D$
with connected $\tilde D$. Let $\tilde\Sigma= \psi^{-1}(\Sigma)$.
Then $\#\tilde\Sigma=2k$  is even.
Now we can construct a  double covering $\phi:D'\to\tilde D$ branched over $\tilde\Sigma$ as  in Lemma \ref{lem:even}.
Then $\Phi=\psi\circ \phi:D'\to D$ is a branched 4-sheet covering with
$$
\chi(D')=2\chi(\tilde D)-2k=2(2\chi(D)-k)<0,
$$
and near each point in $\Sigma'=\Phi^{-1}(\Sigma)$ it is the Levi-Civita transformation.

\medskip

Suppose now that  $D$ is simply connected and $k$ is odd.  Since the case when
$D$ is  a disk was already  covered in Lemma \ref{lem:even}, we may assume that $D=S^2$.
By the condition of Theorem \ref{thm:simple},
$$
k=2m+1>2\chi(S^2)=4,
$$
so $m\ge 2$.

Let $\Sigma=\{a_1,\dots,a_{2m+1}\}$.
By Lemma \ref{lem:even}, there exists a closed surface $N$ and a double covering
$\psi:N\to S^2$ branched over the set $\Lambda=\{a_1,\dots, a_{2m}\}$.
On $N$ the degrees $\alpha_j$ of these singularities
will be replaced by $\alpha'_j=2\alpha_j-2$, so singularities $a'_j=\psi^{-1}(a_j)$, $j=1,\dots,2m$, may disappear,
become weak or remain moderate.
The last  singularity $a_{2m+1}$ will  be replaced by the set $\psi^{-1}(a_{2m+1})=\{b_1,b_2\}$
of two moderate or  Newtonian singularities.
We have
$$
\chi(N)=2\chi(S^2)-2m\le 0,
$$
so the genus of $N$ is $\ge 1$ and there are at least two  moderate or  Newtonian singularities.
This case was already studied. Using Lemma \ref{lem:even} with $\Lambda=\{b_1,b_2\}$, we finally obtain
 a $4$-sheet covering $\phi:X\to S^2$,
branched over $\Sigma$,
such that any point  $q\in \phi^{-1}(\Sigma)$ has multiplicity 2.
Near  $q$ the map $\phi$ is the Levi--Civita regularization.
Now Proposition \ref{prop:reg} is completely proved.

\medskip

Let $\phi:D'\to D$ be the covering in Proposition \ref{prop:reg}.
We have $\#\phi^{-1}(\Sigma)=kK/2$. By the Riemann--Hurwitz formula,
$$
\chi(D')=K\chi(D)-kK/2<0.
$$
Thus we reduced Theorem \ref{thm:simple} to Proposition \ref{prop:weak} also when $k$ is odd.
\qed

\section{Proof of Theorem \ref{thm:top}}

\label{sec:top}

The proof of Theorem \ref{thm:top} is similar to the proof of Theorem \ref{thm:simple}
but instead of the Levi-Civita regularization we use  a sequence of  generalized Levi-Civita regularizations.
We already removed all strong force singularities and assumed that there are no Jacobi singularities. Let
$$
\Sigma=\Delta_{mod}\cup\Delta_{newt}=\cup_{i=1}^m\Delta_{k_i},\qquad \Delta_{k_i}\ne\emptyset, \quad k_i\ge 2.
$$
First we prove Theorem \ref{thm:top} when $D$ is conformally equivalent to a domain in the complex plane.
If  $\Sigma=\Delta_k$ for a single $k$, we can define the regularizing covering as in (\ref{eq:Z}):
\begin{equation}
\label{eq:Z2}
Z=\{(z,w)\in  \mC^2:   w^k=\prod_{a_j\in\Delta_k}(z-a_j)\}.
\end{equation}
Then $\pi:Z\to\mC$ is a  $k$-sheet covering, which is the generalized Levi-Civita regularization of order   
$k$   at each point $a'_j=(a_j,0)$.

In general set
\begin{equation}
\label{eq:X}
Z=\{(z,w_1,\dots,w_m)\in \mC^{m+1}:  w_i^{k_i}=\prod_{a_j\in \Delta_{k_i}}(z-a_j),\; i=1,\dots, m\}.
\end{equation}
It is easy to see that if the points   $a_j$ are distinct, then for $(z,w_1,\dots w_m)\in Z$ the rank of the Jacobi matrix
$$
\frac{\partial(\Phi_1,\dots,\Phi_m)}{\partial(z,w_1,\dots,w_m)},\qquad \Phi_j=w_i^{k_i}-\prod_{a_j\in \Delta_{k_i}}(z-a_j),
$$
equals $m$. Hence $Z$ is a smooth complex curve in $\mC^{m+1}$. The projection
$$
\pi:Z\to \mC,\qquad  \phi(z,w_1,\dots,w_m)=z,
$$
is a covering over $\mC\setminus\Sigma$
with the number of sheets
$$
\#\pi^{-1}(z)=\prod_{i=1}^m k_i=K,\qquad z\notin \Sigma.
$$
The covering is branched over   $\Sigma$: for $a_j\in\Delta_{k_i}$,
$$
\#\pi^{-1}(a_j)=k_1\cdots\hat k_i\cdots k_m=\frac{K}{k_i}
$$
The branching index of every point $q\in \pi^{-1}(a_j)$, $a_j\in \Delta_{k_i}$, equals $\nu(q)=k_i$.
In a neighborhood of $q$ the covering $\pi$ is  the generalized Levi-Civita  transformation (\ref{eq:Levi})
of order $k_i$:
$$
z=a_j+c_j w_i^{k_i}+\dots,\qquad c_j^{-1}=\prod_{a_l\in \Delta_{k_i},\;l\ne j}(a_j-a_l).
$$

Set $X=\pi^{-1}(D)$ and $\phi=\pi|_X$.
We lift the Jacobi metric $g_h$ on $D\setminus\Delta$ to a Riemannian metric  $\phi^*g_h$ on $X\setminus \phi^{-1}(\Delta)$.
The boundary $\partial X=\phi^{-1}(\partial D)$ is geodesically convex.
By Lemma \ref{lem:LC} the order of a singularity  $q\in \pi^{-1}(a_j)$,  $a_j\in\Delta_{k_i}$,
will become $\alpha'_j<\alpha_j$. The singularity
disappears if $\alpha_j=A_{k_i}$ or becomes weak if $A_{k_i}<\alpha_j<A_{k_{i+1}}$.

Let $n_{k_i}=\#\Delta_{k_i}$. By the Riemann--Hurwitz formula,
\begin{eqnarray*}
\chi(X)&=&K\chi(D)-\sum_{i=1}^m\sum_{q\in\pi^{-1}(\Delta_{k_i})} (\nu(q)-1)\\
&=&K\chi(D)-\sum_{i=1}^m n_{k_i}\frac{K}{k_i}(k_i-1)\\
&=&K\chi(D)-\frac K2\sum_{i=1}^m n_{k_i}A_{k_i}= K(\chi(D)-\frac12 A(\Delta)).
\end{eqnarray*}
If (\ref{eq:chi}) holds, then $\chi(X)<0$, so Proposition \ref{prop:weak} works.
This proves Theorem \ref{thm:top} for the case when $D$ is homeomorphic to a plane domain.

\medskip

If $D$ is not a plane domain, we proceed as in the proof of Theorem \ref{thm:simple}.
Take a domain $U\subset D$  containing $\Delta$ and conformally equivalent
to a disk $B(0,1)=\{z:|z|\le 1\}$ in the complex plane.
Define  the covering $\pi:Z\to \mC$ as  in (\ref{eq:X})  and set $Y=\pi^{-1}(U)$.
It is easy to see that the boundary $\partial Y$ is homeomorphic to the curve
$$
S=\{(z,w_1,\dots,w_m)\in \mC^{m+1}: |z|=1,\; w_i^{k_i}=z^{n_{k_i}},\; i=1,\dots, m\}.
$$
Indeed, we can assume that $U=B(0,1)$.  Since $\Sigma\subset U$, for $R\ge 1$ the topology of $Y_R=\pi^{-1}(B(0,R))$
does not depend on $R$. But for large $R$ the boundary $\partial Y_R$ is given by the equations
$w_i^{k_i}=z^{n_{k_i}}(1+O(R^{-1}))$.

If all ratios $n_{k_i}/k_i=d_i$ are integers, then $S$ consists of $K=\prod_{i=1}^m k_i$ closed curves:
$$
S=\cup S_\sigma,\qquad S_\sigma=\{(z,w_1,\dots,w_m):|z|=1,\; w_i=\sigma_i z^{d_i},\; i=1,\dots, m\},
$$
where $\sigma_i$, $\sigma_i^{k_i}=1$, are roots of unity of degree $k_i$ and $\sigma=(\sigma_1,\dots,\sigma_m)$.
The projection $\pi:S_\sigma\to \partial U$ is a diffeomorphism. Hence we can attach to $Y$ a copy of $D\setminus U$
along each circle $S_\sigma$ and obtain a smooth surface $X$ and a branched  $K$-sheet covering $\phi:X\to D$  as above.
We proved:

\begin{lem}\label{lem:X}
Suppose $D$ is a plane domain or all ratios $n_{k_i}/k_i$ are integers. Then there exists an $K=\prod_{i=1}^m k_i$
sheet smooth branched covering $\phi:X\to D$
such that near every point $q\in\phi^{-1}(\Delta_{k_i})$, the map $\phi$ is the generalized Levi-Civita
regularization of order $k_i$. The Euler characteristics of $X$ is given by
$$
\chi(X)=K(\chi(D)-\frac12 A(\Delta)).
$$
\end{lem}

Of course the condition of Lemma \ref{lem:X} is very restrictive. Now consider the general case.
If $D$ is not simply connected,\footnote{Recall that $D$ is oriented, so it is a sphere with handles and holes.}
we first take a covering $\psi:\tilde D\to D$ of degree $K=\prod_{i=1}^m k_i$.
Then $\Sigma$ is replaced by $\tilde\Sigma=\psi^{-1}(\Sigma)$  and  $n_{k_i}$ is replaced by $\tilde n_{k_i}=Kn_{k_i}$.
Then  Lemma \ref{lem:X} works for $D$ and $\Sigma$ replaced by $\tilde D$ and $\tilde\Sigma$.
The composition of coverings $\phi:X\to\tilde D$ and $\psi:\tilde D\to  D$
gives the desired branched covering $\Phi:X\to D$  of degree $K^2$.

At each point in $\phi^{-1}(\Delta_k)$
the covering is the generalized Levi-Civita regularization of order $k$.
Hence all singularities in $\Delta_{reg}$ are regularized, and the remaining singularities are now weak by Lemma \ref{lem:LC}.
The assertion of Theorem \ref{thm:top} follows from Proposition \ref{prop:weak}.

\medskip

If $D$ is simply connected, we may assume that $D=S^2$.  As in Theorem \ref{thm:simple}, this case is  most nontrivial.
The covering $\phi:X\to S^2$ will be different for different types of singularities. In general the lower $A(\Delta)$,
more subtle the construction.

As an example, let us take the lowest possible $A(\Delta)>4$. Then 3 singularities $a_1,a_2,a_3$ are Newtonian,
and the 4th has order $\alpha_4=4/3$, so $A(\Delta)=4\frac13$. Take Newtonian singularities $a_1,a_2$ and perform the
global Levi-Civita regularization in Lemma \ref{lem:even} branched over $\Lambda=\{a_1,a_2\}$.
The Euler characteristics of the new configuration space will be $2\chi(S^2)-2=2$, so it is still a sphere.
Singularities $a_1,a_2$ will disappear, $a_3$ will be replaced by 2 Newtonian singularities $b_1,b_2$,
and  $a_4$ will be replaced by 2 singularities of order $4/3$.

Repeating the regularization in Lemma \ref{lem:even} with   $\Lambda=\{b_1,b_2\}$,
we will have no Newtonian singularities left
but 4 singularities  $c_1,c_2,c_3,c_4$ of order $4/3$. The configuration space is still $S^2$.
Next we perform the regularization $\phi:X\to S^2$  in Lemma \ref{lem:X} of degree $K=3$
branched over singularities $c_1,c_2,c_3$.   The Euler characteristic of the regularized surface $X$ will be
$$
\chi(X)=3(\chi(S^2)-\frac 12 \cdot 3\cdot \frac 43)=0,
$$
so $X=\T^2$. Now singularities $c_1,c_2,c_3$ are regularized, but there are 3 singularities of order $4/3$
in $\phi^{-1}(c_4)$. One more regularization of order 3 branched over these singularities will
get the configuration space $Y$ with  negative Euler characteristics. The composition of all these coverings will
be a branched covering $\Phi:Y\to S^2$ of degree $K=36$. All singularities are now regularized,
and $\chi(Y)<0$.

Other types of singularities on $S^2$ are treated in a similar way. For example, if there are singularities $a_1,a_2,a_3$
with $1\le \alpha<\frac43$ and $a_4$ with $\frac43\le \alpha_4<\frac32$, we take the same 36-sheet covering $\Phi:Y\to S^2$.
On $Y$, all singularities will become weak, so Proposition \ref{prop:weak} applies.

Theorems \ref{thm:top} and \ref{thm:cover} are now proved when there are no Jacobi singularities.

\section{Jacobi singularities}

Until now we assumed that there are no Jacobi singularities.
Suppose we already dealt with other types of singularities  by the generalized Levi-Civita regularization
as described in the previous section  and obtained  a geodesically convex domain $D$  containing only weak and
Jacobi singularities.
Then under the condition of Theorem \ref{thm:top},
$$
\chi(D\setminus \Delta_{jac})=\chi(D)-n<0, \qquad n=\#\Delta_{jac}\ge 1.
$$

For a Jacobi singularity  $a_j$ with $\alpha_j=2$, we have $\rho(x,a_j)=+\infty$ for $x\ne a_j$
as for strong singularities.
Hence the Jacobi metric on $D\setminus \Delta_{jac}$ is complete,
but there is no concavity property of the balls $B_j$ as for strong singularities.
However any nontrivial homotopy class  of closed curves in
 $D\setminus \Delta_{jac}$ contains a minimizer.
The only exception are trivial homotopy classes of small  closed loops $\gamma$ in $B_j$
going $k$ times around a Jacobi singularity $a_j$. The length $J(\gamma)$
of such a loop is bounded: it is close to $2\pi k \sqrt{2m_j}$.
For trivial classes, the infimum of $J$ may be attained on a trivial curve $\gamma\equiv a_j$.
A curve $\gamma$ in a nontrivial homotopy class $\Gamma\subset C^0(S^1,D\setminus\Delta_{jac})$ 
can not be pulled close to a Jacobi singularity $a_j$
without increasing the length $J(\gamma)$ to $+\infty$. Hence $\gamma$ stays in a compact subdomain in
$D\setminus \Delta_{jac}$, so $J$ has a minimum on $\Gamma$. This argument is due to Poincar\'e and Gordon \cite{Gordon}.

So we still get an infinite number of minimal periodic orbits corresponding to nontrivial homotopy classes
in $C^0(S^1,D\setminus \Delta_{jac})$.
Thus the first item of Theorem \ref{thm:top} is proved also in the presence of Jacobi singularities.

To get a compact invariant set with positive topological entropy a different argument is needed,
since we don't have concavity as for strong singularities.
Obtaining chaotic trajectories as limits of minimal periodic geodesics when homotopy classes in $D\setminus\Delta_{jac}$
become more and more complicated will not work since minimal closed geodesics may spiral more and more
close to a Jacobi singularity, so we will not get a compact invariant set.

We will surround pairs of  Jacobi singularities by closed geodesics obtaining a compact geodesically convex set $D'\subset D$.
This requires sufficiently many Jacobi singularities.
We   can achieve this by using   that the order of a Jacobi singularity
does not change under the Levi-Civita regularization, see Lemma \ref{lem:LC}.

If $D$ is not simply connected, then $\chi(D)\le 0$. Take a double covering $\tilde D\to D$ replacing $D$ with
a new domain $\tilde D$ with $\chi(\tilde D)=2\chi(D)$. Now there are $2n\ge 2$ Jacobi singularities
$b_{1},\dots,b_{2n}$.
We take simple curves $C_1,\dots, C_n$ bounding disks  containing only pairs of Jacobi singularities $b_i,b_{n+i}$
and find minimizing  curves $\gamma_i$ in the  homotopy classes $\Gamma_i\subset C^0(S^1,D\setminus\Delta_{jac})$ of curves $C_i$. The minimum on $\Gamma_i$ is attained exist since
curves  $\gamma\in\Gamma_i$ with bounded $J(\gamma)$ can not be pulled  close to $\Delta_{jac}$:
they stay outside small neighborhoods of the singularities.
The minimizers $\gamma_i$ are simple curves and $\gamma_i\cap \gamma_j=\emptyset$, $i\ne j$.
Hence $\gamma_i$ bounds a disk
$D_i$ containing $b_i,b_{i+1}$, and $D_i\cap D_j=\emptyset$ for $i\ne j$.
Replace   $\tilde D$ by a geodesically convex domain $D'=\tilde D\setminus \cup D_i$
with Euler characteristics $\chi(D')=2\chi(D)-n<0$. Now Proposition \ref{prop:weak} works.

If $D$ is simply connected, it is a sphere or a disk. Suppose first that $D$ is a disk with $\chi(D)=1$, then $n\ge 2$.
Perform the Levi-Civita regularization in Lemma \ref{lem:even} for the Jacobi singularities $\Lambda=\{a_1,a_2\}$.
The  domain $D$  will be replaced by $D'$ with $\chi(D')=2\chi(D)-2=0$, so $D'$ is a cylinder $S^1\times [0,1]$.
By Lemma \ref{lem:LC}, the singularities $a_1,a_2$ will be replaced by a pair of  Jacobi singularities $b_1,b_2$.
Taking a minimizing curve around $b_1,b_2$ bounding a disk $U$ as above and deleting  $U$ from $D$ we get again
 a geodesically convex domain $D\setminus U$ with Euler characteristics $\chi(D\setminus U)=\chi(D)-1=-1$.

 Finally, let $D=S^2$, then $n\ge 3$. Perform the Levi-Civita regularization in Lemma \ref{lem:even} for the
 Jacobi singularities $a_1,a_2$.
The configuration space will be replaced by a closed surface $N$ such that $\chi(N)=2\chi(S^2)-2=2$, so $N=S^2$.
But now there are at least 4 Jacobi singularities, two obtained from $a_1,a_2$,
and at least 2 from the remaining one or more singularities.
Using Lemma \ref{lem:even} with 4 singularities, we get new configuration space $\T^2$ and at least 4 Jacobi singularities.
This case was already covered.

Now Theorem \ref{thm:top} is completely proved also in the presence of Jacobi singularities.

\end{document}